\documentclass[a4paper, 12pt]{amsart}
\usepackage{enumerate}
\usepackage{amsmath, amsfonts, amssymb,amsthm}
\usepackage[height=22.5cm, width=15cm]{geometry}
\usepackage[all]{xy}
\usepackage{mathrsfs, yfonts}
\usepackage{graphicx}
\usepackage[T1]{fontenc}
\usepackage[utf8]{inputenc}
\usepackage{hyperref,amssymb,url,upref,verbatim,xspace,mathrsfs}

\newdir^{ (}{{}*!/-5pt/\dir^{(}}
\newdir_{ (}{{}*!/-5pt/\dir_{(}}
\RequirePackage[mathscr]{eucal}
\let\mathcal\mathscr
\newcounter{toto}
\def\thetoto{\arabic{toto}}
\let\oldmarginpar\marginpar
\def\marginpar#1{\refstepcounter{toto}\textsuperscript{\textup{[\thetoto]}}\oldmarginpar{\footnotesize\textsuperscript{[\thetoto]}\,#1}}

\makeatletter
\@namedef{subjclassname@2010}{\textup{2010} Mathematics Subject Classification}
\def\l@section{\@tocline{1}{0pt}{0pc}{}{}}
\def\l@subsection{\@tocline{2}{0pt}{1.5pc}{}{}}
\def\l@subsubsection{\@tocline{3}{0pt}{2pc}{}{}}
\makeatother

\def\shd{\mathcal{D}}
\def\shh{\mathcal{H}}
\def\shj{\mathcal{J}}
\def\shl{\mathcal{L}}
\def\shm{\mathcal{M}}
\def\shn{\mathcal{N}}
\def\sho{\mathcal{O}}
\newcommand{\C}{\mathbb{C}}
\newcommand{\N}{\mathbb{N}}

\newcommand{\Rhom}{R\shh\!om}
\newcommand{\rb}{\mathrm{b}}
\newcommand{\coh}{\mathrm{coh}}

\DeclareMathOperator{\ho}{\mathcal{H}\mathit{om}}

\DeclareMathOperator{\Char}{Char}
\DeclareMathOperator{\coker}{coker}
\DeclareMathOperator{\rD}{\mathsf{D}}

\DeclareMathOperator{\Hom}{Hom}

\let\tilde\widetilde

\let\epsilon\varepsilon

\let\setminus\smallsetminus
\let\leq\leqslant
\let\geq\geqslant

\def\loccit{loc.\kern3pt cit.{}\xspace}
\def\cf{cf.\kern.3em}
\def\eg{e.g.\kern.3em}

\def\resp{\text{resp.}\kern.3em}

\theoremstyle{plain}
\newtheorem{theorem}{Theorem}[section]
\newtheorem{proposition}[theorem]{Proposition}
\newtheorem{lemma}[theorem]{Lemma}

\theoremstyle{definition}
\newtheorem{definition}[theorem]{Definition}

\newtheorem{remark}[theorem]{Remark}

\newtheorem*{claim*}{Claim}

\newcommand{\RedefinitSymbole}[1]{%
\expandafter\let\csname old\string#1\endcsname=#1
\let#1=\relax
\newcommand{#1}{\csname old\string#1\endcsname\,}%
}
\RedefinitSymbole{\forall} \RedefinitSymbole{\exists}

\def\to{\mathchoice{\longrightarrow}{\rightarrow}{\rightarrow}{\rightarrow}}

\def\To#1{\mathchoice{\xrightarrow{\textstyle\kern4pt#1\kern3pt}}{\stackrel{#1}{\longrightarrow}}{}{}}

\let\oldbigoplus\bigoplus
\renewcommand{\bigoplus}{\mathop{\textstyle\oldbigoplus}\displaylimits}
\let\oldbigwedge\bigwedge
\renewcommand{\bigwedge}{\mathop{\textstyle\oldbigwedge}\displaylimits}
\let\oldbigcap\bigcap
\renewcommand{\bigcap}{\mathop{\textstyle\oldbigcap}\displaylimits}
\let\oldbigcup\bigcup
\renewcommand{\bigcup}{\mathop{\textstyle\oldbigcup}\displaylimits}

\newcommand{\parag}[1]{\paragraph{\sc{#1.}} }

\newtheorem{thm}{Theorem}[subsection]

\newtheorem{cor}[thm]{Corollary}

\begin{document}

\setlength{\parindent}{0.em}

\title{On Lisbon integrals}

\date{21/02/20}
\subjclass[2010]{44A99, 32C35, 35A22, 35 A27, 58J15}
\author{Daniel Barlet
\and Teresa Monteiro Fernandes}

\maketitle

\begin{abstract}
We introduce new complex analytic integral transforms, the Lisbon Integrals, which naturally arise in the study of the affine space $\C^k$ of unitary polynomials $P_s(z)$ where $s\in\C^k$ and $z\in \C$, $s_i$ identified to the $i-$th symmetric function of the roots of $P_s(z)$. We completely determine the $\shd$-modules (or systems of partial differential equations) the Lisbon Integrals satisfy and prove that they are their unique global solutions. If we specify a holomorphic function $f$ in the $z$-variable, our construction induces an integral transform which associates a regular  holonomic module quotient of the sub-holonomic module we computed. We illustrate this correspondence in  the case of a $1$-parameter family of exponentials $f_t(z) = exp(t z)$ with $t$ a complex parameter. 

\end{abstract}

\tableofcontents

\newpage

\section{Introduction}
The main purpose of this paper is to understand the behaviour of functions obtained by integration of $\mathscr{C}^{\infty}-$forms on the fibers of a holomorphic proper fibration. This has been investigated by the first author in two extreme cases: when the basis is $1-$dimensional (see \cite{B2}) and when the fibers are finite (see \cite{B3}), but also by many authors in more general settings (see for instance \cite{BMa}, \cite{L}, \cite{S}).\\
In the present paper we look at a simple but very interesting case where the fibers are the roots of the universal monic equation of degree $k$. The general result proved in  \cite{B3} says that the singularity of these functions are controlled by regular holonomic $\shd-$modules. 

Our purpose it to give a precise answer in this special context. For instance, giving two entire functions $f$ and $g$ in $ \mathcal{O}(\C)$ we want to compute the regular holonomic system whose solutions are the $k-$uples of continuous functions on $\C^k$ given by
 $$ \Psi_p(s_1, \dots, s_k) := \sum_{P_s(z_j) = 0} {z_j}^pf(z_j)\bar g(z_j) \quad p \in [0, k-1]$$
where $P_{s}(z) = \sum_{h=0}^k (-1)^hs_hz^{k-h}$ is the universal monic polynomial of degree $k$.\\
As we are interested only in holomorphic derivative in $s$, the function $g$ is irrelevant for the $\shd-$module we are interested in, but, more surprisingly, there exists a sub-holonomic $\shd-$module of which all these $k-$uples of (continuous functions) are solutions ( in the sense of distributions). 

We determine precisely this $\shd-$module, via formula $(@@)$, for which it is useful to consider (in the simplest form, without $g$) the complex integral  representation (\ref{LI2}) of $\Psi_P$.\\
In order to make this computation, a main step is to introduce the trace of differential forms $f(z)ds_1\wedge \dots\wedge ds_{k-1}\wedge dz$ corresponding to the natural holomorphic volume forms on $H := \{ (s, z)\in \C^k\times \C \ / \ P_{s}(z) = 0\}$ identified to $\C^k$ via the map $(s,z) \mapsto (s_1, \dots, s_{k-1}, z)$. These holomorphic traces\footnote{despite the "denominators" in the formula 
$$ Trace(z^pf(z)ds_1 \wedge \dots \wedge ds_{k-1} \wedge dz) = \left(\sum_{P_{s}(z_j)= 0} \frac{{z_j}^pf(z_j)}{P'_{s}(z_j)}\right)ds_1\wedge \dots \wedge ds_k $$ theses forms have no singularity on the discriminant hypersurface $\{\Delta(s) = 0\}$ in $\C^k$.} have a very simple integral representation via the "Lisbon integrals" (see integral representation (\ref{LI1}) below). \\
 Here we explicit the $\shd$ -modules of which (\ref{LI1}) are (the unique) solutions via formula (@) showing that they derive from a very simple one by the usual functorial operations on $\shd$-modules (inverse image and direct image) as follows: \\

Note that the hypersurface $H$ is also defined by the equation $$s_k=(-1)^{k-1}\sum_{h=1}^k(-1)^h s_h z^{k-h}.$$
Let $j: H\subset \C^{k+1}$ denote the closed embedding and let $B_{H|\C^{k+1}}$ denote the regular holonomic $\shd_{\C^{k+1}}$-module of holomorphic distributions supported by $H$. Since the restriction of the projection $$\pi: \C^{k+1}\to \C^k, \, (s,z)\mapsto s$$ to $H$ is proper (with finite fibers),
given a coherent $\shd_{\C^{k+1}}$-module $\shl$ non characteristic for $H$, following \cite{KS4} we obtain a complex in $\rD^\rb_{\coh}(\shd_{\C^k})$, the integral transformed of $\shl$, given by the composition of the usual derived functors of direct image and inverse image for $\shd$-modules, $$D\pi_*(B_{H|\C^{k+1}}\overset{L}{\otimes}_{\sho_{\C^{k+1}}}\shl)\simeq D\pi_*(\shh^1_{[H]}(\shl))\simeq ({\pi|_H})_*j^*\shl$$

Note that $B_{H|\C^{k+1}}=\shh^1_{[H]}(\sho_{\C^{k+1}})$ and that for such a module $\shl$, we have, thanks to Kashiwara's equivalence theorem (cf Th 4.1 and Prop 4.2, \cite{Ka0}),  $\shh^1_{[H]}(\shl)\simeq \shh^0Dj_*Dj^*\shl\simeq j_*j^*\shl$.

We show that the $\shd_{\C^k}$-module determined by all vector functions $\Phi_f$ given by the integral transform $(3)$ ($f$ varying in the space of holomorphic functions in the $z$-variable) is obtained as an integral transform in the sense of Kashiwara and Schapira (\cite{KS4}) of a coherent $\shd_{\C^{k+1}}$-module $\shl$.

 In this note  $\shl$ is the quotient of $\shd_{\C^{k+1}}$ by the ideal generated by the partial derivatives in $s_i$, $i=1,\cdots, k$, hence the sheaf of solutions of $\shl$ is $p^{-1}\sho_{\C}$ where $p(s,z)=z$.
 
 To simplify we shall keep the notation $\pi$ also for the restriction $\pi|_H$ of $\pi$ to $H$.

 As a consequence, we show that Lisbon Integrals (\ref{LI1}) are exactly the global solutions of $\pi_*j^*\shl$.
 
 Moreover, once an entire function $f$ is fixed, we can consider the regular holonomic $\shd_{\C^{k+1}}$-module (denoted by $\shl_f$)  it defines: $$\shl_f=\shd_{\C^{k+1}}/\shj$$ where $\shj$ is the coherent ideal of $\shd_{\C^{k+1}}$ of operators $P$ such that $Pf=0$; hence, according to \cite[Th. 8.1]{Ka3}, $\pi_*j^*\shl_f$ is regular holonomic. We explicit this module in the case of the family $f_t(z)=e^{tz}$ where $t$ is a  complex parameter. \\

\textit{Since integrals (\ref{LI1}) and (\ref{LI2}) are strongly related as explained below, for the sake of simplicity we call both Lisbon Integrals.}\\

We also prove that Lisbon integrals (\ref{LI2}) are global solutions of another $\shd_{
\C^k}$-module $\tilde{\shn}$ which shares with the first this very simple relation: \\
Let $A(s)$ be the $(k, k)-$matrix companion of the unitary polynomial $P_{s}(z)$. Then if $\Phi$ is a solution of $\pi_*j^{*}\shl$ then $\Psi=P'_s(A(s)).\Phi$,  where $P'_s$ denotes the partial derivative of $P_{s}$ with respect to $z$, is a solution of $\tilde{\shn}$. Furthermore, this correspondence $\Phi\leftrightarrow \Psi$ is a bijection when restricting to the  complementary of the discriminant hypersurface  $\{ \Delta(s) = 0 \}$.

Important features of the scalar components of Lisbon integrals are the following:

\begin{itemize}

\item{ They are common solutions of a particular sub-holonomic system. This  aspect will be developed in another paper by the first author. Here we compute only the simplest case $k=2$.}
\item{Each entire function $f$ determines a solution of $\tilde{\shn}$ which scalar component of order $h$ is the trace with respect to $\pi$ in the holomorphic sense of the function $f(z)z^h$ on $H$.}

\end{itemize}

Last but not the least, these computations illustrate the fact that it is not so easy, even in a rather simple situation, to follow explicitly the computations hidden in the ``yoga'' of $\shd$-module theory.

We warmly thank the referee for the many pertinent comments contributing to clarify this work.
\\

\section{ Lisbon Integrals and the differential system they satisfy}

\subsection{Lisbon integrals}

For $(z_{1}, \dots, z_{k}) \in \C^{k}$ denote $s_{1}, \dots, s_{k}$ the elementary symmetric functions of $z_{1}, \dots, z_{k}$. We shall consider in the sequel $s_{1}, \dots, s_{k}$ as coordinates on
$\C^{k} \simeq \C^{k}\big/\mathfrak{S}_{k}$, isomorphism given by the standard symmetric function theorem.\\
We shall denote $P_{s}(z) := \prod_{j=1}^{k} (z - z_{j}) = \sum_{h=0}^{k} (-1)^{h}s_{h}z^{k-h}$ with the convention $s_{0}\equiv 1$.

\textbf{We shall often write  $P(s,z)$ instead of $P_s(z)$ with no risk of ambiguity.}

The discriminant $\Delta(s)$  of $P_{s}$  is the polynomial in $s$ corresponding to the symmetric polynomial  $\prod_{1 \leq i < j \leq k}(z_{i} - z_{j})^{2}$ via the symmetric function theorem.

\begin{lemma}
For $h \in \mathbb{N}$ and $f \in \mathcal{O}(\C)$ any entire holomorphic function, let us define, for $R \gg \vert\vert s\vert\vert$, 
\begin{equation}
 \varphi_{h}(s) := \frac{1}{2i\pi}\int_{\vert \zeta\vert = R} \frac{f(\zeta)\zeta^{h}d\zeta}{P_{s}(\zeta)} .
 \end{equation}
Then $\varphi_{h}(s)$ is independent of the choice of $R$ large enough and defines a holomorphic function on $\C^{k}$. For $\Delta(s) \not= 0$ we have
\begin{equation}
 \varphi_{h}(s) = \sum_{j=1}^{k} \frac{z_{j}^{h}f(z_{j})}{P'_{s}(z_{j})} 
 \end{equation}
where $z_{1}, \dots, z_{k}$ are the roots of $P_{s}(z)$.
\end{lemma}

\parag{Proof} The independence on $R$ large enough when $s$ stays in a compact set of $\C^{n}$ is clear. For $s$ in the interior of a compact set,  $P_{s}(\zeta)$ does not vanish on $\{ \vert \zeta\vert = R \}$ for  $R$ large enough, so we obtain the holomorphy of $ \varphi_{h} $ near any point in $\C^{k}$. The  formula $(2)$ is given by a direct  application of the Residue's theorem.
$\hfill \blacksquare$\\

In fact, it will be convenient to consider the $k$ functions $\varphi_{0}, \dots, \varphi_{k-1}$  as the component of a vector valued function $\Phi := \begin{pmatrix} \varphi_{0}\\ \varphi_{1}\\ \dots \\ \varphi_{k-1}\end{pmatrix}$. Defining  $E(z) := \begin{pmatrix} 1 \\ z \\ \dots \\ z^{k-1} \end{pmatrix}$
 we obtain
\begin{equation}\label{LI1}
 \Phi(s) = \frac{1}{2i\pi}\int_{\vert \zeta\vert = R} \frac{f(\zeta)E(\zeta)d\zeta}{P_{s}(\zeta)} .
 \end{equation}

\begin{definition}\label{Lisbon1}
We call $\Phi$ (sometimes also denoted by $\Phi_f$ when precision is required) the \textit{Lisbon Integral associated to $f$}. The scalar components of $\Phi$, denoted by $\varphi_h$, $h=0,\cdots,k-1$, are called the \textit{scalar Lisbon Integrals}. One also denote by $\varphi_h$ the functions constructed by the same formula, with $h\in\N$, still denominated by "scalar Lisbon Integrals".

\end{definition}

 It will be also interesting to introduce another type of integrals, still called Lisbon Integrals for the sake of simplicity:
 \begin{equation}\label{LI2}
 \Psi(s) :=  \frac{1}{2i\pi}\int_{\vert \zeta\vert = R} \frac{f(\zeta)E(\zeta)P'_{s}(\zeta)d\zeta}{P_{s}(\zeta)}
 \end{equation}
 
 $\Psi$ will also be noted below by  $\Psi_f$ when precision is required.
 
 It is easy to see   that this is again a vector valued holomorphic function on $\C^{k}$ and the Residue's theorem entails that, for $\Delta(s) \not= 0$, the component $\psi_{h}$ of $\Psi$ is given by:
 \begin{equation}
 \psi_{h}(s) = \sum_{j=1}^{k} z_{j}^{h}f(z_{j}) .
 \end{equation}

\begin{proposition}\label{inj} If $f$ is not identically zero then $\Phi$ and $\Psi$ are non zero vector-valued holomorphic functions on $\C^k$.
 \end{proposition}
 
 \parag{Proof} Suppose that $f$ is non  identically zero. Then the statement follows as an immediate consequence of the non vanishing of the Van der Monde determinant of $z_{1}, \dots, z_{k}$ when these complex numbers are pair-wise distinct.$\hfill \blacksquare$
 
  \parag{An example} Take $f \equiv 1$. Then formula (5) shows that $\psi_{h}(s)$ is the $h-$th Newton symmetric functions of the roots of the polynomial $P_{s}$. So it is a quasi-homogeneous polynomial in $s$ of  weight $h$ (the weight of $s_{j}$ is $j$ by definition).\\
    Let us show that we have $\varphi_{h}(s) \equiv 0$ for $h \in [0, k-2]$ and $\varphi_{k-1}(s) \equiv 1$ in this case. For $h \in [0, k-2]$ the formula (1) gives the estimate (with $f \equiv 1$)
 $$ \vert \varphi_{h}\vert \leq  \frac{R^{h+1}}{(R - a)^{k-1}} $$
 if each root of $P_{s}$ is in the disc $\{ \vert z \vert \leq a\}$  when  $R > a > 0$. When $R \to +\infty$ this gives $\varphi_{h}(s) \equiv 0$ for $h \in [0, k-2]$.\\
 For $h = k-1$ write
  $$kz^{k-1} = P'_{s}(z) - \sum_{h=1}^{k-1} (-1)^{h}(k-h)s_{h}z^{k-h-1}.$$
  This gives, using the previous case and formula (2), that $\varphi_{k-1}(s) \equiv 1$.\\

 \subsection{The partial differential system}

 Let us introduce the $(k, k)$ matrix $A$ (the companion matrix) associated to the polynomial $P_{s}$:
 
 \bigskip
 
\begin{equation}\label{EA}
 A := \begin{pmatrix} 0 & 1 & 0 & \cdots & \quad  & 0 \\ 0 & 0 & 1& 0 & \cdots & 0 \\ \quad  \\ 0 & \quad & \cdots & \quad &\quad & 0 \\ \quad \\ \quad  \\ 0 & \quad & \cdots & \quad & 0 & 1    \\ (-1)^{k-1}s_{k}& \cdots & (-1)^{h-1}s_{h}& \cdots &\cdots& s_{1} \end{pmatrix}
 \end{equation}

\begin{thm}\label{system}
{\it The vector valued holomorphic function $\Phi$ on $\C^{k}$ satisfies the following differential system}
\begin{equation*}
 (-1)^{k+h}\frac{\partial \Phi}{\partial s_{h}}(s) = \frac{\partial (A^{k-h}\Phi)}{\partial s_{k}}(s) \qquad \forall s \in \C^{k} \quad {\rm and} \quad \forall h \in [1, k-1].\tag{@}
 \end{equation*}
{\it Moreover, this system is integrable\footnote{We shall explain in the proof what we mean here.} and if $\Phi$ is a solution of this system, so is $A\Phi$}.
\end{thm}

The proof of this result will use several lemmas.

\begin{lemma}\label{elem.1}
Let $A$ be a $(k, k)$ matrix with entries in $\C[x]$ and put $B := \lambda\frac{\partial A}{\partial x}$ where $\lambda$ is a complex number. Let $M$ be the $(2k, 2k)$ matrix given by
$$ M := \begin{pmatrix} A & B \\ 0 & A \end{pmatrix} .$$
Then for each $p \in \mathbb{N}$ we have 
\begin{equation*}
 M^{p} = \begin{pmatrix} A^{p} & B_{p} \\ 0 & A^{p}\end{pmatrix} \tag{a}
 \end{equation*}
where $B_{p} := \lambda\frac{\partial (A^{p})}{\partial x} $.
\end{lemma}

\parag{Proof} As the relation $({a})$ is clear for $p = 0, 1$ let us assume that it has been proved for $p$ and let us prove it for $p+1$. We have:
$$ \begin{pmatrix} A^{p} & B_{p}\\ 0 & A^{p}\end{pmatrix}\begin{pmatrix} A & B \\ 0 & A \end{pmatrix} = \begin{pmatrix} A^{p+1} & A^{p}B + B_{p}A \\ 0 & A^{p+1} \end{pmatrix} $$
which allows to conclude.$\hfill \blacksquare$\\

\begin{cor}\label{simple 1}
{\it For each integer  $p \in [0, k-1]$ the following equality holds in the module $\C^{k}\otimes_{\C} \C[s_{1}, \dots, s_{k}, z]\big/(P^{2})$ over the $\C-$algebra $\C[s_{1}, \dots, s_{k}, z]\big/(P^{2})$} 

\begin{equation}\label{E7}
 z^{p}E(z) = A^{p}E(z) + (-1)^{k-1}P_{s}(z)\frac{\partial (A^{p})}{\partial s_{k}}E(z)
 \end{equation}
{\it In particular, for any entire function $f$ (of the variable $z$), we have $$\Phi_{zf}=A(s)\Phi_f$$ }
 
{\it Moreover the following identity in the module $\C^{k}\otimes_{\C} \C[s_{1}, \dots, s_{k}, z]\big/(P^{2})$  holds}
\begin{equation}\label{E7-bis}
P'_{s}(z)E(z) = P'_{s}(A)E(z) +  (-1)^{k-1}P_{s}(z)\frac{\partial (P'_{s}(A))}{\partial s_{k}}E(z). \
\end{equation}

\end{cor}

\parag{Proof} In the basis $1, z, \cdots, z^{k-1}, P_{s}(z), zP_{s}(z), \cdots, z^{k-1}P_{s}(z)$ of this algebra which is a free rank $2k$ module on  $\C[s_{1}, \dots, s_{k}]$, the multiplication by $z$ is given by the matrix $M$ of the previous lemma with  $A$ as in (\ref{EA}) and with $B := (-1)^{k-1}\frac{\partial A}{\partial s_{k}}$. This proves equality (\ref{E7}).\\
As $P'_{s}(z) = \sum_{h=0}^{k-1} (-1)^{h}(k-h)z^{k-h-1} $ does not depend on $s_{k}$ it is enough to sum up the previous equalities with $p = k-h-1$ with the convenient coefficients to obtain the equality (\ref{E7-bis}).
$ \hfill \blacksquare$\\

  
 
\begin{lemma}\label{simple 2} 
For any $h \in [1, k]$ and any $p \in \mathbb{N}$ the matrix $A$ in $(6)$ satisfies the relation:
\begin{equation}\label{E2}
(-1)^{k-h} \frac{\partial A^{p}}{\partial s_{h}} = \frac{\partial A^{p}}{\partial s_{k}}A^{k-h}
\end{equation}
\end{lemma}

\parag{Proof} The case $p=1$ of $(\ref{E2})$  is an easy direct computation on the matrix $A$. Assume that the assertion is proved for $p \geq 1$. Then Leibnitz's rule gives:
\begin{align*}
& (-1)^{k-h}\frac{\partial A^{p+1}}{\partial s_{h}} = (-1)^{k-h} \frac{\partial A^{p}}{\partial s_{h}}A + A^{p}(-1)^{k-h}\frac{\partial A}{\partial s_{h}} \\
& \quad =  \frac{\partial A^{p}}{\partial s_{k}}A^{k-h+1} + A^{p}\frac{\partial A}{\partial s_{k}}A^{k-h} = \frac{\partial A^{p+1}}{\partial s_{k}}A^{k-h} \
\end{align*}
concluding the proof of (\ref{E2}).$\hfill \blacksquare$\\

\parag{Proof of the theorem \ref{system}} By derivation inside the integral in  (\ref{LI1}) we obtain:
\begin{align*}
& \frac{\partial \Phi}{\partial s_{h}}(s) =  \frac{1}{2i\pi} \int_{\vert \zeta\vert = R} \ f(\zeta)E(\zeta)(-1)^{h+1}\zeta^{k-h}\frac{d\zeta}{P_{s}(\zeta)^{2}} \quad {\rm and \ in \ particular} \\
&  \frac{\partial \Phi}{\partial s_{k}}(s) =  \frac{1}{2i\pi} \int_{\vert \zeta\vert = R} \ f(\zeta)E(\zeta)(-1)^{k+1}\frac{d\zeta}{P_{s}(\zeta)^{2}} 
\end{align*}
Now for $h \in [1, k-1]$ we use the formula of corollary \ref{simple 1} to obtain:
$$ \frac{\partial \Phi}{\partial s_{h}} = (-1)^{h+1}A^{k-h}(-1)^{k-1} \frac{\partial \Phi}{\partial s_{k}} + (-1)^{h+1}(-1)^{k-1}\frac{\partial A^{k-h}}{\partial s_{k}}\Phi $$
that is to say, we obtain $(@)$ as desired:
\begin{equation*}
(-1)^{k+h}\frac{\partial \Phi}{\partial s_{h}} =  \frac{\partial (A^{k-h}\Phi)}{\partial s_{k}} \qquad \forall  h \in [1, k] 
\end{equation*}

By the integrability of the system $(@)$ we mean that for any $\Phi$  such that $(@)$ holds, then the computation of the partial derivatives $\frac{\partial^{2}\Phi}{\partial s_{h}\partial s_{j}}$ using the system $(@)$ gives  a symmetric result  in $(h, j)$ for any pair $(h, j)$ in $[1, k]$. Note that if $h$ or $j$ is equal to $k$ the assertion is trivial.\\
 So consider a couple $(h, j) \in [1, k-1]^{2}$. Thanks to Lemma \ref{simple 2} we have :

\begin{align*}
& (-1)^{h+j}\frac{\partial^{2} \Phi}{\partial s_{j}\partial s_{h}} = (-1)^{k-j}\frac{\partial }{\partial s_{k}}\big[\frac{\partial (A^{k-h}\Phi) }{\partial s_{j}}\big] \\
& \qquad =  (-1)^{k-j}\frac{\partial }{\partial s_{k}}\big[\frac{\partial A^{k-h}}{\partial s_{j}}\Phi + A^{k-h}\frac{\partial \Phi}{\partial s_{j}} \big] \\
& \qquad = \frac{\partial }{\partial s_{k}}\big[ \frac{\partial A^{k-h}}{\partial s_{k}}A^{k-j}\Phi + A^{k-h}\frac{\partial (A^{k-j}\Phi)}{\partial s_{k}}\big] \\
& \qquad =  \frac{\partial^{2} }{\partial s_{k}^{2}}\big[A^{2k-h-j}\Phi\big] 
\end{align*}
which is symmetric in $(h, j)$.\\
To finish the proof of the theorem we have to show that $A\Phi$ is a solution of $(@)$ when $\Phi$ is a solution of $(@)$. This is given by the following computation 

\begin{align*}
& (-1)^{k-h}\frac{\partial (A\Phi)}{\partial s_{h}} =   (-1)^{k-h}\frac{\partial A}{\partial s_{h}}\Phi +  (-1)^{k-h}A\frac{\partial \Phi}{\partial s_{h}} \\
& \qquad = \frac{\partial A}{\partial s_{k}}.A^{k-h}\Phi + A\frac{\partial A^{k-h}\Phi}{\partial s_{k}} = \frac{\partial (A^{k-h}(A\Phi))}{\partial s_{k}}
\end{align*}
which also uses Lemma \ref{simple 2}.$\hfill \blacksquare$\\

\parag{Remark} A consequence of our computation on the integrability of the system $(@)$ is the fact that for any solution $\Phi$ and any pair $(h, j) \in [1, k]$ the second order partial derivative $\frac{\partial^{2} \Phi}{\partial s_{j}\partial s_{h}}$ only depends on $h+j$. This implies that any scalar  Lisbon integral $\varphi_{h} $ satisfies
\begin{equation}\label{E3}
\frac{\partial^{2} \varphi_{h}}{\partial s_{p}\partial s_{q+1}} = \frac{\partial^{2} \varphi_{h}}{\partial s_{p+1}\partial s_{q}} \qquad \forall \ p, q \ {\rm such \ that} \    1 \leq p < q \leq k-1 
\end{equation}

 Let us denote by $\Delta$ the discriminant hypersurface  $\Delta=\{ \Delta(s) = 0 \}$. An easy calculation shows that away of $\Delta$ the matrix $P'_s(A(s))$ is invertible.
 
 The  next corollary of theorem \ref{system}  gives an analogous system to $(@)$  for the vector function $\Psi$ defined in $(4)$ which is singular along $\Delta$.

\begin{cor}\label{system 2}

\begin{enumerate}
\item {The vector valued holomorphic function $\Psi$ on $\C^{k}$ satisfies the following differential system:
\begin{equation*}
(-1)^{k+h}\frac{\partial \Psi}{\partial s_{h}}(s) = \frac{\partial (A^{k-h}\Psi)}{\partial s_{k}}(s) + (-1)^{k}(k-h)A^{k-h-1}P'_{s}(A)^{-1}\Psi(s)   \tag{@@}
\end{equation*}
$\forall h \in [1, k-1]$,  which is singular along the discriminant hyperdurface $$ \Delta := \{ s \in \C^{k}; \Delta(s) = 0 \}.$$}
\item{If $\Psi$ is any solution of $(@@)$ then $A\Psi$ is also a solution of $(@@)$.}
\item{If $\Phi$ is any solution of $(@)$ then $\Psi=P'_s(A)\Phi$ is a solution of $(@@)$.}
\item{If $\Psi$ is any solution of $(@@)$ on $\C^k \setminus \Delta$ then $\Phi:=P'_s(A(s))^{-1}\Psi$ is a solution of $(@)$ on $\C^k \setminus \Delta$. }
\end{enumerate}

\end{cor}
Statement $(1)$ follows by (\ref{E7-bis}) in Corollary \ref{simple 1}.

Statement $(3)$ follows by direct computation: \\
For such a $\Psi=P'_s(A)\Phi$, for each $h \in [1, k-1]$
\begin{equation*}
 (-1)^{k-h}\frac{\partial \Psi}{\partial s_{h}}(s) = (-1)^{k-2h}(k-h)A^{k-h-1}\Phi(s) + \sum_{p=0}^{k-1} (-1)^{p}(k-p)(-1)^{k-2h}s_{p}\frac{\partial (A^{k-p-1}\Phi)}{\partial s_{h}}(s)
 \end{equation*}
 and using now the fact that $A^{k-p-1}\Phi$ is solution of $(@)$ we obtain
 \begin{align*}
 &  (-1)^{k-h}\frac{\partial \Psi}{\partial s_{h}}(s) = (-1)^{k-2h}(k-h)A^{k-h-1}\Phi(s) + \sum_{p=0}^{k-1} (-1)^{p}(k-p)s_{p}\frac{\partial (A^{2k-p-h-1}\Phi)}{\partial s_{k}}(s)\\
 & =  (-1)^{k-2h}(k-h)A^{k-h-1}\Phi(s) + \frac{\partial (A^{k-h}P'_{s}(A)\Phi)}{\partial s_{k}}(s) \\
 & =  (-1)^{k-2h}(k-h)A^{k-h-1}P'_{s}(A)^{-1}\Psi(s) + \frac{\partial (A^{k-h}\Psi)}{\partial s_{k}}(s)
 \end{align*}
 which gives the formula $(@@)$.

Since $P'_{s}(A)$ commutes with $A$, the  assertion $(2)$ is easy.$\hfill \blacksquare$\\
The proof of $(4)$ is  a simple consequence of $(\ref{E2})$ in Lemma \ref{simple 2}.

\subsection{Example: The case $k = 2$}

In this example we will explicit the system $(@)$ and also a partial differential operators in the Weyl algebra $\C[s_{1}, s_{2}]\langle \partial_{s_{1}}, \partial_{s_{2}}\rangle$, which annihilates  the scalar components of its solutions. The left ideals in $\shd_{\C^k}$  of differential operators annihilating respectively  the scalar components of the solutions of $(@)$ and of $(@@)$  for arbitrary $k$ are described in \cite{B1}.\\ 
 \vspace{2mm}
Here we use the notations $s := s_{1} $ and $p := s_{2}.$
In that case, with $\Phi=(\varphi_0, \varphi_1)$, the differential system (@) becomes:

\begin{equation}\label{E11}
 \frac{\partial \varphi_{0}}{\partial s} = \frac{\partial \varphi_{1}}{\partial p} 
 \end{equation}
 \begin{equation}\label{E12}
\frac{\partial \varphi_{1}}{\partial s} = - \varphi_{0} -p\frac{\partial \varphi_{0}}{\partial p} + s\frac{\partial \varphi_{1}}{\partial p}
\end{equation}
corresponding to the matrix  $A := \begin{pmatrix} 0 & 1 \\ -p & s \end{pmatrix}$. Differentiating  (\ref{E12}) with respect to  $p$ we obtain, after substituting via (\ref{E11})
\begin{align*}
&   \frac{\partial^{2} \varphi_{0}}{\partial s^{2}} = -2 \frac{\partial \varphi_{0}}{\partial p}  - p \frac{\partial^{2} \varphi_{0}}{\partial p^{2}}   - s \frac{\partial^{2} \varphi_{0}}{\partial s \partial p} \\
& {\rm and \ so} \qquad  \frac{\partial^{2} \varphi_{0}}{\partial s^{2}} + s \frac{\partial^{2} \varphi_{0}}{\partial s \partial p} + p\frac{\partial^{2} \varphi_{0}}{\partial p^{2}}  + 2 \frac{\partial \varphi_{0}}{\partial p} = 0 \tag{$\sharp$}
\end{align*}
Differentiating (\ref{E12}) with respect to  $s$ we obtain, after substituting via (\ref{E11})
\begin{align*}
& -  \frac{\partial^{2} \varphi_{1}}{\partial s^{2}} =  \frac{\partial \varphi_{1}}{\partial p} + p \frac{\partial^{2} \varphi_{1}}{\partial p^{2}} +  \frac{\partial \varphi_{1}}{\partial p} + s\frac{\partial^{2} \varphi_{1}}{\partial s \partial p} \\
&   {\rm and \ so} \qquad    \frac{\partial^{2} \varphi_{1}}{\partial s^{2}} + s \frac{\partial^{2} \varphi_{1}}{\partial s \partial p} + p \frac{\partial^{2} \varphi_{1}}{\partial p^{2}}  + 2 \frac{\partial \varphi_{1}}{\partial p} = 0 \tag{$\sharp\sharp$}
\end{align*}
Then the second order  differential operator of weight $-2$  
\begin{equation*}
\Theta := \frac{\partial^{2}}{\partial s^{2}} + s\frac{\partial^{2}}{\partial s\partial p} + p\frac{\partial^{2}}{\partial p^{2}} + 2\frac{\partial}{\partial p}   \tag{$\sharp\sharp\sharp$}
\end{equation*}
anihilates  $\varphi_{0}$ and $\varphi_{1}$ for any solution $\Phi$ of the system ((\ref{E11}), (\ref{E12})). \\

\parag{A direct proof that $\Theta$ anihilates scalar Lisbon Integrals for all $m \in \mathbb{N}$} We have, for any entire function $f : \C \to \C$ and for $R \gg \max\{ \vert s\vert, \vert p\vert\}$:
\begin{equation*}
\varphi_{m}(s, p) = \frac{1}{2i\pi}\int_{\vert \zeta\vert = R} \ f(\zeta)\frac{\zeta^{m}d\zeta}{\zeta^{2}-s\zeta + p} \tag{a}
\end{equation*}

This gives:
\begin{align*}
& \frac{\partial \varphi_{m}}{\partial s}(s, p) = \frac{1}{2i\pi}\int_{\vert \zeta\vert = R} \ f(\zeta)\frac{\zeta^{m+1}d\zeta}{(\zeta^{2}-s\zeta + p)^{2}}  \tag{b} \\
& \frac{\partial \varphi_{m}}{\partial p}(s, p) = - \frac{1}{2i\pi}\int_{\vert \zeta\vert = R} \ f(\zeta)\frac{\zeta^{m}d\zeta}{(\zeta^{2}-s\zeta + p)^{2}}  \tag{c} \\
& \frac{\partial^{2} \varphi_{m}}{\partial s^{2}}(s, p) = 2\frac{1}{2i\pi}\int_{\vert \zeta\vert = R} \ f(\zeta)\frac{\zeta^{m+2}d\zeta}{(\zeta^{2}-s\zeta + p)^{3}}  \tag{d} \\
& \frac{\partial^{2} \varphi_{m}}{\partial s\partial p}(s, p) = -2\frac{1}{2i\pi}\int_{\vert \zeta\vert = R} \ f(\zeta)\frac{\zeta^{m+1}d\zeta}{(\zeta^{2}-s\zeta + p)^{3}} \tag{e} \\
&  \frac{\partial^{2} \varphi_{m}}{\partial p^{2}}(s, p) = 2\frac{1}{2i\pi}\int_{\vert \zeta\vert = R} \ f(\zeta)\frac{\zeta^{m}d\zeta}{(\zeta^{2}-s\zeta + p)^{3}} \tag{f}
\end{align*}
Now it is easy to check that $ (d) + s(e) + p(f) + 2(c) = 0$.\\

\section{The left action of $\Gamma(\C, \shd_{\C})$ on Lisbon integrals}
For each entire function $f$ of the variable $z$ we shall henceforward denote by $\Phi_f$ the associated Lisbon integral (previously generically denoted by $\Phi$). Clearly the assignement $$f\mapsto \Phi_f$$ is $\C$-linear and, according to Proposition \ref{inj}, it is injective.
\begin{lemma}\label {laction} 

Let $f$ be an entire function on $\C$ and let $\Phi_{f}$ the corresponding Lisbon integral.
\begin{enumerate}
\item{Let $g$ be an entire function of $z$. Then $\Phi_{gf}=g(A(s))\Phi_f$. In particular $\Phi_f=f(A(s))\Phi_1$}.

\item{We have the identity
$$ \Phi_{\partial_{z}(f)}(s) = -\nabla\Phi_{f}(s) + \big(\sum_{h=0}^{k-1} (k-h)s_{h}\partial_{s_{h+1}}  \big)(\Phi_{f})(s) $$
where $\nabla$ is the constant $(k, k)$ matrix $\begin{pmatrix} 0 & 0 & 0 & \dots & 0 \\ 1 & 0 & 0 & \dots & 0 \\ 0 & 2 & 0 & \dots & 0 \\ 0 & 0 & & \dots  &  0 \\ 0 & \dots &0 & k-1 & 0 \end{pmatrix} $.}
\end{enumerate}

\end{lemma}

\parag{Proof}  

When $g$ is a polynomial on $z$,  statement $(1)$ follows easily in Corollary \ref{simple 1}.
For an arbitrary entire function $g$, it is a consequence of \cite[Lem 3.1.8]{BaM}.

Let us prove $(2)$:
Consider the Lisbon integral
$$ \Phi_{\partial_{z}f}(s) =  \frac{1}{2i\pi}\int_{\vert \zeta\vert = R} f'(\zeta)E(\zeta)\frac{d\zeta}{P_{s}(\zeta)} $$
Integration  by parts gives, as $\partial_{z}(E)(z) = \nabla E(z)$ :
\begin{equation*}
\Phi_{\partial_{z}(f)}(s) = - \nabla\Phi_{f}(s)  +  \frac{1}{2i\pi}\int_{\vert \zeta\vert = R} f(\zeta)E(\zeta)\frac{P'_{s}(\zeta)}{P_{s}(\zeta)^{2}}d\zeta \tag{*}
\end{equation*}
Now, using the equalities $P'_{s}(\zeta) = \sum_{h=0}^{k-1} (-1)^{h}(k-h)s_{h}\zeta^{k-h-1}$ and
$$ \frac{\partial \Phi_{f}}{\partial s_{h}}(s) = - \frac{1}{2i\pi}\int_{\vert \zeta\vert = R} f(\zeta)E(\zeta)\frac{(-1)^{h}\zeta^{k-h}d\zeta}{P_{s}(\zeta)^{2}} $$
we obtain the formula of the lemma.$\hfill \blacksquare$\\

\parag{Remark} From formula $(^{*})$ and according to $(\ref{E7-bis})$
we obtain also the formula
\begin{equation*}
\Phi_{\partial_{z}(f)}(s) = - \nabla\Phi_{f}(s)  + (-1)^{k-1}\frac{\partial (P'_{s}(A)\Phi_{f})}{\partial s_{k}}(s) \tag{**}
\end{equation*}

\smallskip
Note that the formula $\Phi_{(z.f)'} = \Phi_{f} + \Phi_{zf'}$ corresponding to the usual relation $\partial_{z}z - z\partial_{z} = 1$ follows from the linearity of the map $f\mapsto \Phi_f$ and the Leibniz rule $(zf)'=f+zf'$.

It is not obvious that when $\Phi$ is solution of the system $(@)$, then
 $$(\Phi \partial_{z})(s) := -\nabla\Phi(s)  + (-1)^{k-1}\partial_{s_k}(P'_{s}(A)\Phi)(s)$$
 (given by the formula $(^{**})$) is also a solution of the same system. A direct verification of this fact is consequence of the formula given in the lemma below.

\begin{lemma}\label{verif}
We have the following identity:
$$  \nabla A^p - A^p\nabla + pA^{p-1} = (-1)^{k-1}\partial_{s_k}(A^p)P'_{s}(A), \quad \forall p \geq 1, \  \forall s \in \C^{k}.$$
\end{lemma}

\parag{Proof} Let check the case $p = 1$ first. It is an easy computation to obtain that $ \nabla A - A\nabla + Id $ is the matrix which have all lines equal to $0$ excepted its last one which is given by $(x_{1}, \dots, x_{k})$ with 
$ x_{h} = (-1)^{k-h}hs_{k-h}$ for $h \in [1, k]$ with $s_{0} \equiv 0$. On the other hand, the matrix $(-1)^{k-1}\partial_{s_k}A$ has only a non zero term at the place $(k, 1)$ which equal to $1$,
so it is quite easy to see that $(-1)^{k-1}(\partial_{s_k}A)A^{p}$ has only a non zero term at the place $(k, p+1)$ with value $(-1)^{k-1}$. According to the computation of $(-1)^{k-1}\partial_{s_k}(A)P'_{s}(A)$
we conclude the  desired formula for $p = 1$\\
Now we shall make an induction on $ p \geq 1$ to prove the general case. Then assume $p \geq 2$ and  the formula proved for $p-1$. Then write
$$ \nabla A^p - A^p \nabla = (\nabla A^{p-1} - A^{p-1} \nabla) A + A^{p-1}(\nabla A - A\nabla) .$$
Using the induction hypothesis and the case $p = 1$ gives
\begin{align*}
& \nabla A^p - A^p\nabla = -(p-1) A^{p-2} A - A^{p-1} +\\
&  \qquad  + (-1)^{k-1}\partial_{s_k}(A^{p-1}) P'_s(A)A + (-1)^{k-1}A^{p-1}\partial_{s_k}(A)P'_s(A) \\
& = -pA^{p-1} + (-1)^{k-1}(\partial_{s_k}(A^{p-1})A + A^{p-1}\partial_{s_k}(A))P'_s(A) \\
&= -p A^{p-1} + (-1)^{k-1}\partial_{s_k}(A^p)P'_s(A).\hfill \qquad \qquad \qquad \qquad \qquad   \blacksquare\\
\end{align*}

Now we shall make the direct verification that $\Phi$ solution of $(@)$ implies that 
 $$-\nabla \Phi + (-1)^{k-1}\partial_{s_k }(P'_{s}(A)\Phi)$$ 
 is also solution of $(@)$ :
\begin{align*}
& X :=  (-1)^{k+h}\partial_{s_h}(\nabla\Phi) - \partial_{s_k}(A^{k-h}\nabla\Phi) = \partial_{s_k}( \nabla A^{k-h}\Phi) - \partial_{s_k}(A^{k-h}\nabla \Phi) \\
& = \partial_{s_k}\big[-(k-h)A^{k-h-1}\Phi + (-1)^{k-1}\partial_{s_k}(A^{k-h})P'_s(A)\Phi \big] 
\end{align*}
thanks to the previous lemma. Also, using the fact that $P'_s(A)\Phi$ is a simple linear combination of solutions of $(@)$ (with non constant coefficients, but very simple), we obtain the formula :
$$ (-1)^{k+h}\partial_{s_h}(P'_s(A)\Phi) = (-1)^k(k-h)A^{k-h-1}\Phi + \partial_{s_k}(A^{k-h}P'_s(A)\Phi) .$$

Then:

\begin{align*}
& Y :=  (-1)^{k+h}\partial_{s_h}(\partial_{s_k}(P'_s(A)\Phi)) - \partial_{s_k}(A^{k-h}\partial_{s_k}(P'_s(A)\Phi)) \\
& = \partial_{s_k}\big[(-1)^k(k-h)A^{k-h-1}\Phi + \partial_{s_k}(A^{k-h}P'_s(A))\Phi) - A^{k-h}\partial_{s_k}(P'_s(A)\Phi)\big] \\
& = \partial_{s_k}\big[(-1)^k(k-h)A^{k-h-1}\Phi + \partial_{s_k}(A^{k-h})P'_s(A)\Phi)\big]
\end{align*}

and $-X + (-1)^{k-1}Y = 0$, as desired.$\hfill \blacksquare$\\

\section{The $\shd_{\C^k}$-module associated to Lisbon Integrals}

We shall begin by recalling some basic facts on  $\shd$-modules and by fixing notations.\\

For a morphism of manifolds $f:Y\to X$, we use the notation of \cite{KS1} $$f_d:= {^t}f': T^*X\times_X\times Y\to T^*Y$$ and
 $$f_{\pi}:T^*X\times_X\times Y\to T^*X$$
  the associated canonical morphisms of vector bundles.

We recall that a conic involutive submanifold $V$ of the cotangent bundle $T^*Z$ of a manifold $Z$ (real or complex) is $\mathit{regular}$ if the restriction $\omega|_V$ of the canonical $1$-form $\omega$ on $T^*Z$ never vanishes outside the $0-$section. Recall also that if $(x_1,\cdots,x_n, \xi_1,\cdots,\xi_n)$ are canonical symplectic coordinates on $T^*Z$, then $\omega(x,\xi)=\sum_{i=1}^n\xi_idx_i$.
  
Let us fix some $k\in\N$, $k\geq 2$. In $\C^{k+1}=\C^k\times \C$ we consider the coordinates $(s_1,\cdots, s_k, z)$. As in the previous sections we set $$P(s,z)=z^k+\sum_{h=1}^{k-1}(-1)^hs_hz^{k-h}$$ Obviously $$P(s_1,\cdots, s_k, z)=0\Longleftrightarrow s_k=(-1)^{k-1}\sum_{h=0}^{k-1}(-1)^hs_hz^{k-h},$$  where $s_0=1$.
We note $s=(s_1,\cdots, s_k)$ and $s':=(s_1,\cdots,s_{k-1})$.

Let $H$ be the smooth hypersurface of $\C^{k+1}$ given by the zeros of $P(s, z)$. 

Let us denote by $\shl$ the $\shd_{\C^{k+1}}$-module with one generator $u$ defined by the equations $\partial u/\partial{s_1}=\cdots=\partial u/ \partial{s_k}=0$. 
Such module is an example of a so called \textit{partial de Rham systems}, which have the feature, among others, that their characteristic varieties are non singular regular involutive. In our case we have
$$\Char\shl = \{(s, z); (\eta, \tau)\in \C^{k+1}\times \C^{k+1} \text{such that} \, \eta = 0 \} .$$
Since $H\subset\C^{k+1}$ is defined by the equation 
$$P(s, z)=(-1)^ks_k+\sum_{h=0}^{k-1}(-1)^hs_hz^{k-h}=0,$$
 $T^*_{H}\C^{k+1}$ is the subbundle of $T^*\C^{k+1}|_H$ described by
  $$\{(s,z;\eta, \tau), (s,z)\in H, \exists \lambda \in \C \quad {\rm such \ that} \quad  (\eta, \tau)=\lambda dP(s,z)\}.$$
  This means that $P(s, z) = 0$ and that their exists $\lambda \in \C$ with $ \eta_{h} = \lambda(-1)^{h}z^{k-h}$ for each $h \in [1, k]$ and that $\tau = \lambda P'(s, z)$.
  Hence as  $\eta_{k} = \lambda(-1)^{k}$ and this implies:
  $$\Char\shl\cap T^*_{H}\C^{k+1}\subset T^*_{\C^{k+1}}\C^{k+1}$$
   (as usual, $T^*_{\C^{k+1}}\C^{k+1}$ denotes the zero section of $T^*\C^{k+1}$), in other words $H$ is non characteristic for $\shl$. Let us denote by $j:H\hookrightarrow \C^{k+1}$ the closed embedding. By Kashiwara's classical results (which can be found in \cite{Ka2}) it follows that  the induced system $D j^*\shl$ by $\shl$ on $H$ is concentrated in degree zero and $\shn:=\shh^0 Dj^*\shl$ is a  $\shd_{\C_H}$-coherent module whose characteristic variety is exactly
$$j_dj_{\pi}^{-1}\Char(\shl).$$

Recall that
 $$D j^*\shl:=\sho_H\overset{L}{\otimes}_{j^{-1}\sho_{\C^{k+1}}}j^{-1}\shl$$ and in this non-characteristic case 
 $$\simeq j^{-1}\frac{\sho_{\C^{k+1}}}{P(s',s_k,z)\sho_{\C^{k+1}}}\otimes_{j^{-1}\sho_{\C^{k+1}}}j^{-1}\shl$$

We have
 $$\shn:= j^{-1}(\frac{\shd_{\C^{k+1}}}{P\shd_{\C^{k+1}}+\shd_{\C^{k+1}}\partial_{s_1}+\cdots+\shd_{\C^{k+1}}\partial_{s_k}})\simeq j^{-1}(\sho_{\C^{k+1}}/P\sho_{\C^{k+1}})<\partial_z>$$ which is isomorphic as a $\shd_H$-module to $$\sho_H<\partial_z> \simeq \frac{\shd_H}{\shd_H\partial_{s_1}+\cdots+\shd_H\partial s_{k-1}}$$ where $\partial_{s_i}$ stands for the derivation $\partial/\partial s_i$ on $\sho_H$ and $\partial_z$ as a derivation on $\sho_H$ is the class of $\partial_z$ in the quotient above.

 In particular  $\shn$ is sub-holonomic and it is a partial de Rham system similarly to $\shl$.

Let us now determine the image of $\shn$ under the morphism $\pi: \C^k\simeq H\to \C^k$ given by $(s',z)\mapsto (s', s_k)$.
Clearly $\pi$ is proper surjective with finite fibers.

Recall that one denotes by $\shd_{\C^k\leftarrow H}$ the transfer $(\pi^{-1}\shd_{\C^k}, \shd_H)$-bimodule
 $$(\pi^{-1}\shd_{\C^{k}}\otimes_{\pi^{-1}\sho_{\C^{k}}}\pi^{-1}\Omega_{\C^{k}}^{\otimes^{-1}})\otimes_{\pi^{-1}\sho_{\C^{k}}}\Omega_H$$

Recall also that, according to the properness and the fiber finiteness of $\pi$, we have 

$$D\pi_*\shn\simeq \shh^0D\pi_*\shn=\pi_*(\shd_{\C^k\leftarrow H}\otimes_{\shd_H}\shn)$$
where we abusively use the notation $\pi_*$ for the direct image functor in the categoy of $\shd$-modules in the two left terms and for the direct image functor for sheaves in the right term.
According to  \cite[Th. 4.25 and 4.27]{Ka2} (see also the comments in loc.cit. before Theorem 4.27), one knows that $D\pi_*\shn$ is concentrated in degree zero and that 
$$\Char\shh^0 D\pi_*\shn=\pi_{\pi}\pi_d^{-1}\Char \shn.$$

So we may henceforward denote for short $\pi_*\shn:=D\pi_*\shn$ without ambiguity.

Let $\Delta$ as above be the zero set of the discriminant of $P$, which can also be defined as the image by $\pi$ of the subset of $H$ defined by $\{ P'(s, z) = 0\}$.\\
Since $\pi_d$ is given by the $k\times k$ matrix 
$$\begin{pmatrix}
      Id &  0  \\
     \partial s_k/\partial {s'}&\partial s_k/\partial z\\ 
\end{pmatrix}^{T}$$ 
we conclude that $\Char \pi_*\shn$ is the image by $\pi_{\pi}$ of the set
$$ \{ (s', z); (\eta', \tau) \in \C^{k}\times \C^{k} \ / \  \eta_{j} = -(-z)^{k-j-1}\tau, \quad \forall j \in [1, k-1] \}$$
so it is given by the set
$$ \{(s, \eta) \in \C^{k}\times \C^{k} \ / \ \exists z \in \C \ {\rm such \ that} \ P_{s}(z)= 0 \ {\rm and \ with} \ \eta_{j} = (-z)^{k-j-1}\eta_{k} \quad \forall j \in [1, k-1] \}$$

Then  $\Char \pi_*\shn$  is an involutif analytic subset of $T^*\C^k$ with codimension $ k-1$ which proves the following:

\begin{lemma}\label{subhol}$\pi_*\shn$ is a subholonomic $\shd_{\C^{k+1}}$-module.

\end{lemma}
\begin{remark}
Let $\tilde{\shn}$ denote the $\shd_{\C^{k}}$-module associated to $(@@)$. Then $\tilde{\shn}$ is clearly not subholonomic.
\end{remark}
 

\begin{proposition}\label{Basic}
{\it The $\shd_{\C^k}$-module $\pi_{*}\mathcal{N}$ is the quotient of $\shd_{\C^k}^{k} \simeq \shd_{\C^k}\otimes_{\C}\C^{k}$ by the action of}
$$ \mathcal{A}_{h} := \partial_{s_{h}} \otimes Id_{\C^{k}} + (-1)^{k-h-1}\partial_{s_{k}}\otimes A(s)^{k-h} \quad {\rm for} \ j \in [1, k-1] .$$
{\it Moreover the action of $z$ and $\partial_{z}$ on $\pi_{*}\mathcal{N}$ deduced from the action of $\shd_H$ on $\mathcal{N}$\footnote{Note that $z$ and $\partial_{z}$ commute with $\partial_{s_{h}}$ for $h \in [1, k-1]$ in $\shd_H$.} are given respectively by}
$$ \mathcal{A}_{0} := 1 \otimes A(s) \quad {\rm and} \quad \mathcal{B} := 1\otimes \nabla + (-1)^{k-1}\partial_{s_{k}} \otimes P'_{s}(A(s)) $$
{\it where we put $P'_{s}(z) := (\partial_{z}(P_{s}(z))$ and }
$ \nabla :=  \begin{pmatrix} 0 & 0 & \dots & 0\\ 1 & 0 & \dots & 0 \\ 0 & 2 & 0 & \dots \\ \dots & \dots & \dots & \dots \\ 0 & \dots & k-1 & 0 \end{pmatrix} $
\end{proposition}

\parag{Proof}  Our goal is to explicit $\pi_*\shn$ and to check that it coincides with the  $\shd_{\C^k}$-module associated to the  system (@) in Theorem \ref{system}.

In a first step we explicit the transfer-module $$\shd_{\C^k\leftarrow H}:=\pi^{-1}\shd_{\C^{k}}\otimes_{\pi^{-1}\sho_{\C^{k}}}(\pi^{-1}\Omega_{\C^{k}}^{\otimes^{-1}}\otimes_{\pi^{-1}\sho_{\C^{k}}}\Omega_H)$$ as a $(\pi^{-1}\shd_{\C^k}, \shd_H)$-bimodule.\\
 The next step is to determine the cokernel of $$\alpha: (\shd_{\C^k\leftarrow H})^{k-1}\to \shd_{\C^k\leftarrow H}$$ 
$$(u_1,\cdots,u_{k-1})\mapsto \sum_{i=1}^{k-1}u_i\partial_{s_i}$$



The last step is to apply $\pi_*$.

Let us denote for short $$\sigma:=\omega_H\otimes \omega_{\C^k}^{\otimes^{-1}}:=ds_1\wedge\cdots\wedge ds_{k-1}\wedge d_z\otimes (ds_1\wedge\cdots\wedge ds_{k-1}\wedge ds_k)^{\otimes^{-1}}$$ 
the generator of the line bundle $ \pi^{-1}\Omega_{\C^{k}}^{\otimes^{-1}}\otimes_{\pi^{-1}\sho_{\C^{k}}}\Omega_H$.\\

Recall that $\sho_H=\sho_{\C^{k+1}}/\shj$, where $\shj$ is the ideal generated by $P(s,z)$. Hence $\sho_H$ is a $\pi^{-1}\sho_{\C^k}$-free module with rank $k$ since each section $a(s',z)$ of $\sho_H$ is equivalent, by Weierstrass Division Theorem, to a unique polynomial $\,\sum_{j=0}^{k-1}a_j(s',s_k)z^j$ modulo $P(s,z)$, for some sections $a_j$ of $\sho_{\C^k}$.

Hence $\shd_{\C^k\leftarrow H}$ is a left $\pi^{-1}\shd_{\C^k}$-free module of rank $k$ generated by the $k-$sections \\ $(1\otimes z^j\sigma)_{j=0,\cdots,k-1}$. Since the right action of each operator in $\shd_H$ is $\pi^{-1}\shd_{\C^{k}}$-linear, it is sufficient to calculate each $(1\otimes z^j\sigma)\partial _{s_i}, i=1,\cdots,k-1, \,j=0,\cdots, k-1$.

Now recall that $H$ is defined in $\C^{k+1}$ by the equation $s_k=(-1)^{k-1}\sum_{h=0}^{k-1}(-1)^hs_hz^{k-h}$ with the convention $s_{0} = 1$  and so $s_{1}, \dots, s_{k-1}, z$ are global coordinates on $H$. Then we have in $H$
\begin{equation*}
\frac{\partial s_{k}}{\partial s_{h}} = (-1)^{k-h-1}z^{k-h} \quad {\rm and } \quad \frac{\partial s_{k}}{\partial z} = (-1)^{k-1}P'_{s}(z) 
\end{equation*}
where $P'_{s}(z)$ does not depend on $s_{k}$.\\
 Let $F := 1 \otimes E(z)\sigma $ denote the basis $(1\otimes z^{j}\sigma), j \in [0, k-1]$ of the free rank $k$ left $\pi^{-1}(\shd_{\C^{k}})-$module $\shd_{\C^{k}\leftarrow H}$.
Recall that, according to \cite[Rem. 4.18]{Ka2}, in view of the generators described above, the action of $\shd_H$ in $\shd_{\C^k\leftarrow H}$  is defined by the following formulas, where we consider $F$ as a  $k-$vector and use the usual matrix product
\begin{align*}
&  F\theta(s') = \theta(s')F \quad {\rm where} \quad \theta \in \mathcal{O}_{H} \quad {\rm does \ not \ depend \ on} \ z    \tag{0*} \\
&  Fz = \   A(s)F    \tag{1*}  \\
& - F\partial_{s_{h}} = \partial_{s_{h}}F + (-1)^{k-h-1}\partial_{s_{k}}(A(s)^{k-h}F) \quad \forall h \in [1, k-1]   \tag{2*} \\
& - F\partial_{z} =  \nabla F + (-1)^{k-1}\partial_{s_{k}}\big(P'_{s}(A(s))F\big)  \tag{3*}
\end{align*}
where we have used the equalities $zE(z) = A(s)E(z)$ and $\partial_{z}(E(z)) = \nabla E(z)$.\\

Summing up:
\begin{itemize}
\item{For $g\in\sho_H$ represented by \  $\sum_{r=0}^{k-1}g_r(s)z^r,$ \ the $(k\times k)$ matrix $G$ of the \\
 $\pi^{-1}\shd_{\C^k}$-linear morphism defined by $g$ on $\shd_{\C^k\leftarrow H}$ is given by \\
 $G := \sum_{r=0}^{k-1} g_{r}(s)A(s)^{r}$}
\item{Let us consider the $\shd_{\C^k}$-linear morphism $\alpha : (\shd_{\C^k}^{k})^{(k-1)}\to \shd_{\C^k}^k$ defined by the following  $k-1$ $(k, k)-$matrices}
\end{itemize}
$$ \mathcal{A}_{h} := \partial_{s_{h}}\otimes Id_{\C^{k}} + (-1)^{k-h-1}\partial_{s_{k}}\otimes A(s)^{k-h} $$

Let  $\Phi$ be in $\mathcal{O}_{\C^{k}}^{k}$. In view of the relation $(2^ *)$, the map $(1\otimes z^{h}\sigma) \mapsto \Phi_{h}$, for $ h \in [0, k-1]$ will induce an element of $\ho_{D_{\C^{k}}}(\pi_{*}\mathcal{N}, \mathcal{O}_{\C^{k}})$, that is to say a solution of $\pi_{*}\mathcal{N}$, if and only if we have $\partial_{s_{h}}(\Phi) = (-1)^{k-h}\partial_{s_{k}}(A^{k-h}\Phi)$, that is, in and only if $\Phi$ satisfies (@), since the generator of $\pi_{*}\mathcal{N}$ is anihilated by the action of  $\partial_{s_{h}}$ for each $h \in [1, k-1]$.\\
In conclusion:
$$ \pi_{*}\mathcal{N} \simeq \coker \alpha$$
 by the finitness of the fibers of $\pi$.
 $\hfill \blacksquare$\\



\remark\label{R}
$\shn$ is naturally endowed with a structure of right $\Gamma(\C;\shd_{\C})-$module. By functoriality $\pi_*\shn$ is also a right $\Gamma(\C;\shd_{\C})-$module and its structure coincides with the induced by the right $\shd_{\C}$ action defined by $(1^*)$ and $(3^*)$ on  
$\shd_{\C^{k}\leftarrow H}\otimes_{\shd_H}\shn$, since it commutes with each $\partial_{s_i}$, for $i=1,\cdots,k-1$. Therefore we obtain a natural left action of $\Gamma(\C, \shd_{\C})$ on  $\ho_{\shd_{\C^k}}(\pi_*\shn, \sho_{\C^k})$.

We also conclude, according to Lemma \ref{laction}:
\begin{proposition}\label{P:D}
The left action of $\Gamma(\C; \shd_{\C})$ defined by the above Remark \ref{R} on $\ho_{\shd_{\C^k}}(\pi_*\shn, \sho_{\C^k})$ coincides with the left action of $\Gamma(\C; \shd_{\C})$ on Lisbon integrals. 
\end{proposition}


\bigskip

If $\Phi$ is a solution of $\pi_*\shn$, replacing in the formula $(3^*)$ the second term thanks to the equality obtained for $\Phi$ after applying $(2^*)$ or equivalently (@), we also derive a right action of $\partial_{z}$ which is given by the formula
$$-\Phi\partial_{z} = \nabla\Phi - \sum_{h=0}^{k-1} (k-h)s_{h}\partial_{s_{h+1}}\Phi $$
 
\bigskip

Our next goal is to conclude in Proposition \ref{P4} below that there are no global holomorphic solutions of $\pi_*\shn$ other than those of the form $\Phi_f$, for some holomorphic function $f$ only depending on $z$.
Since $j$ is non-characteristic we have an isomorphism
 $$\,j^{-1}\Rhom_{\shd_{\C^{k+1}}}(\shl, \sho_{\C^{k+1}})\simeq \Rhom_{\shd_{H}}(\shn,\sho_{H})$$

According to Theorem 4.33 (2)  of \cite{Ka2}, making $X=H, Y=\C^k, f=\pi,  \shn=\sho_{\C^k}$ in loc.cit, we obtain 

\begin{theorem}\label{T:1}
For any open subset $\Omega$ of $\C^k$ we have an isomorphism functorial on $\shn$ compatible with restrictions to open subsets$$R\Gamma (\pi^{-1}(\Omega); \Rhom_{\shd_{H}}(\shn, \sho_{H}))\simeq R\Gamma (\Omega; \Rhom_{\shd_{\C^{k}}}(\pi_*\shn,\sho_{\C^{k}}))$$
\end{theorem}

Recall that this isomorphism uses as a tool the "trace morphism": $\pi_*\sho_H\to \sho_{\C^k}$ constructed in \cite[Prop. 4.34]{ Ka2}.

Since for any open subset $\Omega$ and any $\shd_H$-module $\mathcal{P}$, $\Gamma(\Omega; \cdot)$ and $\ho_{\shd_H}(\mathcal{P}, \cdot)$ are left exact functors, since if $\Omega$ is a Stein open set and if $\mathcal{P}$ admits a global resolution by free $\shd_H$-modules of finite rank, then $\Rhom_{\shd_H}(\mathcal{P}, \sho)$ is represented by a complex in degrees $\geq 0$ with $\Gamma(\Omega,\cdot)$-acyclic entries,
we conclude that
$$H^0(R\Gamma(H; \Rhom_{\shd_{H}}(\shn, \sho_{H})))=\Gamma(H; \ho_{\shd_{H}}(\shn, \sho_{H})) \quad {\rm and} $$
$$H^0(R\Gamma(\C^k; \Rhom_{\shd_{\C^k}}(\pi_*\shn, \sho_{\C^k})))=\Gamma(\C^k; \ho_{\shd_{\C^k}}(\pi_*\shn, \sho_{\C^k}))$$ therefore Theorem \ref{T:1} entails
a $\C$-linear isomorphism
$$ T: \Hom_{\shd_{H}}(\shn, \sho_{H})\simeq\Hom_{\shd_{\C^k}}(\pi_*\shn, \sho_{\C^k})$$

\begin{proposition}\label{P4}
The correspondence $$f\mapsto \Phi(f):=\Phi_f$$ 

defines a $\C$-linear isomorphism $$\Phi: \Gamma(\C;\sho_{\C})\to \Gamma(\C^k; \ho_{\shd_{\C^k}}(\pi_*\shn, \sho_{\C^k}))=\Hom_{\shd_{\C^k}}(\pi_*\shn, \sho_{\C^k})$$ 
Moreover, this isomorphism are also $\Gamma(\C, \shd_{\C}$)-left linear.
\end{proposition}
\parag{Proof}

The last statement is clear from Proposition \ref{P:D} and Lemma \ref{laction}.

The remaining of the statement is equivalent to prove that $\Phi: f\mapsto \Phi_f$ defines an isomorphism $\Hom_{\shd_{H}}(\shn,\sho_{H})\to \Hom_{\shd_{\C^k}}(\pi_*\shn,\sho_{\C^k})$.

We already know that $\Phi$ is injective. It remains to prove that $\Phi$ is surjective.
For each $f\in\Gamma(\C, \sho_{\C})$, we introduce the regular holonomic $\shd_{\C}$-module $\shm_f$ (a regular flat holomorphic connection on $\C$) of which the constant sheaf $\C f$ in degree zero is the complex of holomorphic solutions. 

Note that $\shn\simeq \sho_{\C^{k-1}}\boxtimes \shd_{\C}$ where we consider $\C^{k-1}$ endowed with the coordinates $(s_1,\cdots,s_{k-1})$ and $\C$ with the coordinate $z$.
 We denote by $\shn_f$ the regular holonomic $\shd_H$-module (a regular flat holomorphic connection on $H$)

$$\shn_f:=\sho_{\C^{k-1}}\boxtimes \shm_f.$$

It is clear that $\shn_f$ is a quotient of $\shn$,  and, by the left exactness of $\pi_*$,  $\pi_*\shn_f$ is a quotient of $\pi_*\shn$. Moreover, according to Proposition \ref{P:D} and Lemma \ref{laction}, 
$\Phi_f$ belongs to $\Hom_{\shd_{\C^k}}(\pi_*\shn_f, \sho_{\C^k})$. According to Theorem \ref{T:1}, for each $f$ we have a $\C$-linear isomorphism $T_f:\Hom_{\shd_H}(\shn_f, \sho_H)\simeq \Hom_{\shd_{\C^k}}(\pi_*\shn_f, \sho_{\C^k})$ which makes this last one a one dimensional $\C$-vector space. Moreover, by left exactness of $\Hom$ and the exactness of $\pi_*$, we have monomorphisms $\Hom_{\shd_{\C^k}}(\pi_*\shn_f, \sho_{\C^k}){\subset} \Hom_{\shd_{\C^k}}(\pi_*\shn, \sho_{\C^k})$ and, by functoriality, we have $T(f)=T_f(f)$.

We shall use the following result:
\begin{lemma}\label{P2}
Suppose that $f\neq 0$. Then
 $\Hom_{\shd_{\C^k}}(\pi_*\shn_f,\sho_{\C^{k}})$ is a one dimensional $\C$-vector space generated by $\Phi_f$. 
\end{lemma}

\parag{Proof}
The result follows by Proposition \ref{inj} since $\Phi_f$ is a non zero element of $\Hom_{\shd_{\C^k}}(\pi_*\shn_f,\sho_{\C^{k}})$ hence it is a generator as a $\C$-vector space.
$\hfill \blacksquare$\\

Let us now end the proof of Proposition \ref{P4}.\\
Clearly $\Hom_{\shd_H} (\shn,  \sho_H)=\sum_{f}\Hom_{\shd_H}(\shn_f, \sho_H)$ and, according to Lemma \ref{P2}, for each $f$, $\Hom_{\shd_{\C^k}}(\pi_*\shn_f,\sho_{\C^{k}})$
is the $\C$-vector space spanned by $\Phi_f$;  hence $T(f)=\lambda \Phi_f$ for some $\lambda\in\C^*$. Since $\Phi_{\lambda f}=\lambda \Phi_f$ we conclude that $\Phi$ is surjective which gives the desired result. 
$\hfill \blacksquare$\\

As a consequence, isomorphism $\Phi$ explicits isomorphism of Theorem \ref{T:1} since they coincide up to the multiplication by a constant $\lambda\neq 0$.

\subsection{An example}To conclude this article, let us give an interesting example of choice of the entire function $f$ on $\C$ for which we explicit the regular holonomic system on $\C^{k}$ associated to the corresponding Lisbon integrals.\\

\parag{The case $f_{t}(z) := e^{tz}$} Let us fix a parameter $t \in \C^{*}$ and consider the entire function $f_{t}(z) := e^{t z}$.\\
 First remark that according to Lemma \ref{laction}, we have $ \frac{\partial E(z)}{\partial z} = \nabla E(z)$ where $\nabla$ is the $(k, k)$ matrix given by
$$ \nabla :=  \begin{pmatrix} 0 & 0 & \dots & 0\\ 1 & 0 & \dots & 0 \\ 0 & 2 & 0 & \dots \\ \dots & \dots & \dots & \dots \\ 0 & \dots & k-1 & 0 \end{pmatrix} $$
We have \footnote{Remember that $t$ is a fixed complex parameter}
\begin{equation*}
 \Phi_{f_t}(s) := \frac{1}{2i\pi}\int_{\vert \zeta\vert = R} \frac{e^{t\zeta}E(\zeta)d\zeta}{P_{s}(\zeta)}  
 \end{equation*}
and, according to the linearity of $\Phi_{(\cdot)}$, we also have $$t\Phi_{f_t}(s)=\Phi_{\partial_z(f_t)}(s)$$ which, applying $(^{\ast\ast})$, entails

\begin{equation}\label{E4}
 t\Phi_{f_t}(s)
= -\nabla\Phi_{f_t}(s) + (-1)^{k-1}\frac{\partial (P'_{s}(A)\Phi_{f_t})}{\partial s_{k}}(s)
\end{equation}
Hence
\begin{equation}\label{E6}
(t Id + \nabla)\Phi_{f_t}(s) = (-1)^{k-1}\frac{\partial (P'_{s}(A)\Phi_{f_t})}{\partial s_{k}}(s).
\end{equation}

This also implies the following  equation for $\Psi_{f_t}$ away of the discriminant hypersurface  $\Delta$:
\begin{equation}
  (tId + \nabla)P'_{s}(A)^{-1}\Psi_{f_t}(s) = (-1)^{k-1}\frac{\partial (\Psi_{f_t})}{\partial s_{k}}(s)
  \end{equation}
  for
  $$ \Psi_{f_t}(s) = P'_{s}(A) \Phi_{f_t}(s) = \frac{1}{2i\pi}\int_{\vert \zeta\vert = R} \frac{e^{t\zeta}P'_{s}(\zeta)E(\zeta)d\zeta}{P_{s}(\zeta)} .$$
   
Combining $(\ref{E6})$ with the system $(@)$  it is easy to see that we obtain a meromorphic integrable connexion on the trivial vector bundle of rank $k$ on $\C^{k}$ with a pole along the discriminant hypersurface. 

The regularity of this meromorphic connexion is then consequence of the regularity of the $\shd_{H}-$module $\mathcal{N}_{e^{tz}}=\shd_{H}u$ which is given by the equations 
 $$\partial_{s_{h}}u = 0, \quad \forall  h \in [1, k-1] \quad  {\rm  and} \quad  (\partial_{z} - t)u= 0$$
 which is clearly regular holonomic on $H$. So its direct image by $\pi$ (as a $\shd_{\C^{k}}-$module) is regular holonomic on $\C^{k}$  (see \cite[Th. 8.1]{Ka3}).

\bigskip

\thanks{The research of T. Monteiro Fernandes was supported by Funda\c c\~ao para a Ci\^encia e a Tecnologia, UID/MAT/04561/2019.}\\

\address{Daniel Barlet,\\
 Institut Elie Cartan, G\'eom\'{e}trie, Universit\'e de Lorraine, Nancy, France,\\
 CNRS UMR 7502   and  Institut Universitaire de France}\\
email: \textit{daniel.barlet@univ-lorraine.fr}\\

\address{Teresa Monteiro Fernandes\\ Centro de Matem\'atica e Aplica\c c\~oes Fundamentais-CIO \\ and \\ 
Departamento de Matem\' atica da Faculdade de Ci\^encias da Universidade de Lisboa,\\ 
Bloco C6, Piso 2, Campo Grande, 1749-016, Lisboa  Portugal}\\
email : \textit{mtfernandes@fc.ul.pt}

\end{document}